\documentclass{amsart}

\usepackage{eucal}
\input{xypic}
\xyoption{all}

\begin{document}

\newcommand{\C}{{\mathbb C}}
\newcommand{\Q}{{\mathbb Q}}
\newcommand{\R}{{\mathbb R}}
\newcommand{\T}{{\mathbb T}}
\newcommand{\Z}{{\mathbb Z}}

\newcommand{\Diff}{{\rm Diff}}
\newcommand{\e}{{\bf e}}
\newcommand{\Emb}{{\rm Emb}}
\newcommand{\gr}{{\rm Grass}}
\newcommand{\GG}{{\check \mathbb G}}
\newcommand{\half}{{\textstyle{\frac{1}{2}}}}
\newcommand{\hx}{{\surd x}}
\newcommand{\M}{{\overline{\mathcal M}}}
\newcommand{\res}{{\rm res}}
\newcommand{\tr}{{\rm Trace}}

\title {Tate cohomology of circle actions as a Heisenberg group} 
\author{Jack Morava}
\address{Department of Mathematics, Johns Hopkins University, Baltimore,
Maryland 21218}
\email{jack@math.jhu.edu}
\thanks{The author was supported in part by the NSF}
\date {15 September 2001}
\begin{abstract} 
We study the Madsen-Tillmann spectrum $\C P^\infty_{-1}$ as a quotient 
of the Mahowald pro-object $\C P^{\infty}_{-\infty}$, which is closely 
related to the Tate cohomology of circle actions. That theory has an
associated symplectic structure, whose symmetries define the Virasoro 
operations on the cohomology of moduli space constructed by Kontsevich 
and Witten. \end{abstract}

\maketitle
\section {Tate cohomology of circle actions} \bigskip

\noindent
{\bf 1.1} If $E$ is a geometric bordism theory (such as integral homology),
its Tate cohomology $t^*_\T E$ can be constructed by tom Dieck stabilization
from the geometric theory of $E$-manifolds with $\T$-action, the action 
required to be free on the boundary [8]. If $E$ is multiplicative, so is 
$t^*_\T E$; there is a cofibration sequence
\[
\cdots \to E^*B\T_+ \to t^*_\T E \to E_{-*-2}B\T_+ \to \cdots \;
\]
in which the boundary map sends a $\T$-manifold with boundary to the
quotient of its boundary by the (free) $\T$-action. When $E$ is 
complex-oriented [eg $MU$ or $H\Z$] this sequence reduces to
a short exact sequence which identifies $t^*_\T E$ with the Laurent 
series ring $E^*((\e))$ obtained by inverting the Euler class in $E^*(B\T) 
= E^*[[\e]]$, and the boundary map can be calculated as a formal residue; 
more precisely, the formal Laurent series $f$ maps to the residue
of $f d \log_E$ at $\e = 0$, where $d \log_E$ is the invariant differential 
of the formal group law of $E$. When $E$ is not complex-orientable, $t_\T E$ 
can behave very differently [4], as the Segal conjecture shows: up to a 
profinite completion,
\[
t_\T S^0 \sim S^0 \vee S^1 \prod B \T/C
\]
where $C$ runs through proper subgroups of $\T$. 
\bigskip

\noindent
There is a related but simpler theory $\tau^*_\T E$ defined by 
manifolds with free $\T$-action on the boundary alone, which fits in an exact 
sequence
\[
\cdots \to E^* \to \tau^*_\T E \to E_{-*-2}B\T_+ \to \cdots;
\]
ignoring the interior $\T$-action defines a truncation
\[
t^*_\T E \to \tau^*_\T E \;.
\]
{\bf 1.2} It is natural to think of $t^*_\T E$ as the $E$-homology of 
a version [1] of Mahowald's pro-spectrum $\C P^\infty_{-\infty}$, constructed
from the inverse system 
\[
\{ \C P^\infty_{-k} = Th(-k\eta) \}
\]
of Thom spectra defined by the filtered vector bundle
\[
\cdots \subset k \eta \subset (k+1) \eta \subset \cdots
\]
over $\C P^\infty$. To be more precise,  
\[
t_\T E \sim E \wedge S^2 \C P^\infty_{-\infty}
\]
as spectra, and there is a similar equivalence
\[
\tau_\T E \sim E \wedge S^2 \C P^\infty_{-1} \;.
\]
From this point of view, the morphism from $t^*_\T E$ to $\tau^*_T E$
is the $E$-homology of the collapse map 
\[
\C P^\infty_{-\infty} \to \C P^\infty_{-1} \;.
\]
The cobordism class defined by a family of complex-oriented surfaces with 
a free circle action on its boundary, parametrized by $X$, defines an 
element of 
\[
\tau^{-2}_\T MU(X_+) = [X_+, \C P^\infty_{-1} \wedge MU] \;;
\]
the image of this class under the Thom homomorphism from $MU$ to $H\Z$ 
is the homomorphism
\[
H^*(\C P^\infty_{-1},\Z) \to H^*(X_+,\Z)
\]
defined by the classifying map of Madsen and Tillmann, sketched below.
\bigskip

\noindent
{\bf 1.3} When $E = H\Z$, the symmetric bilinear form 
\[
f,g \mapsto (f,g) = \res_{\e = 0} (f g \; d\e)
\]
on the Laurent series ring $t_\T H\Z$ is nondegenerate, and the involution 
$I(z) = z^{-1}$ on $\T$ defines a symplectic form
\[
\{f,g\} = (I(f),g)
\]
which restricts to zero on the subspaces of elements of degree $\geq 0$ 
and $< 0$. The Tate cohomology thus has an intrinsic inner product, with 
canonical polarization and involution. \bigskip

\noindent
The functor $\GG$ which sends the commutative ring $A$ to the
set of formal Laurent series 
\[
\GG(A) = \{ g = \sum_{k \gg 0} g_k x^{k+1}  \in A((x)) \:|\: g_0 \in 
A^{\times} \;,\; g_k \in \surd A \; {\rm if} \; 0 > k \}
\]
(i.e. with $g_0$ a unit and $g_k$ nilpotent when $k$ is negative) in fact 
takes values in the category of groups, with formal composition of series
as the operation. This group of invertible `nil-Laurent' series has a 
linear representation on the abelian group valued functor $A \mapsto A((x))$, 
but it is a little too large to be conveniently representable; in a certain 
sense it is an ind-pro-algebraic analog of the group of diffeomorphisms 
of the circle. \bigskip

\noindent
{\bf 1.4} It is tempting to interpret $\GG$ as a group of automorphisms of the
Tate cohomology, but the most obvious action does not preserve the
symplectic structure. Kontsevich-Witten theory suggests a better alternative: 
there is an embedding   
\[
\e^k \mapsto \gamma_{-k-\half}(x) \;:\; t^*_\T H\Z \to \R((\hx))
\]
of symplectic modules, defined using the fractional divided power 
\[
\gamma_s(x) = \frac {x^s}{\Gamma(1+s)} \;,
\]
in which the symplectic structure on the target is defined by
\[
u,v \mapsto \{u,v\} = \res_{x=0} \; u d v \;. 
\]
[The reals are a notational convenience: some powers of $\pi$ have been
ignored.] Over any field of characteristic zero, the square root of an 
invertible nil-Laurent series in $x$ is an invertible {\it odd} nil-Laurent 
series in $\hx$, and it makes better sense to think of the group $\GG^-_{1/2}$
of such series [10 \S 1.3] as symplectic automorphisms of $t^*_\T H\Q \subset 
\R((\hx))$. The half-integral shift comes ultimately from the fact that 
$t_\T H\Z$ is not Spanier-Whitehead self-dual; rather, its dual is most 
naturally interpreted as its own double suspension. 

\section{Madsen-Tillmann and Kontsevich-Witten}

\noindent
Madsen and Tillmann construct a map
\[
\coprod_{g \geq 0} B\Diff(F_g) \to \Omega^\infty \C P^\infty_{-1}
\]
which is compatible with gluing of surfaces; in particular, it defines a 
lax functor from the two-dimensional topological gravity category [11] to 
a topological category with one object and an $H$-space of morphisms. 
[Reversing the orientation of a surface corresponds to the involution $I$.]
The point of this note is to identify a suitable subgroup of $\GG^-_{1/2}$ 
as the motivic automorphisms of this functor. \bigskip

\noindent
{\bf 2.1} Here is a quick account of one component of [7]: if $F \subset 
\R^n$ is a closed two-manifold embedded smoothly in a high-dimensional 
Euclidean space, its Pontrjagin-Thom construction $\R^n_+ \to F^{\nu}$
maps compactified Euclidean space to the Thom space of the normal bundle 
of the embedding. The tangent plane to $F$ is classified by a map
$\tau : F \to \gr_{2,n}$to the Grassmannian of oriented two-planes in 
$\R^n$, and the canonical two-plane bundle $\eta$ over this space has a 
complementary $(n-2)$-plane bundle, which I will call $(n - \eta)$. The 
normal bundle $\nu$ is the pullback along $\tau$ of $(n - \eta)$; composing 
the map induced on Thom spaces with the collapse defines 
\[
\R^n_+ \to F^\nu \to \gr_{2,n}^{(n - \eta)} \;.
\]
The space $\Emb(F)$ of embeddings of $F$ in $\R^n$ becomes 
highly connected as $n$ increases, and the group $\Diff(F)$ of
orientation-preserving diffeomorphisms of $F$ acts freely
on it, defining a compatible family
\[
\R^n_+ \wedge_{\Diff} \Emb(F) \to \gr_{2,n}^{(n - \eta)}
\] 
which can be interpreted as a morphism
\[
B\Diff(F) \to \lim \; \Omega^n \gr_{2,n}^{(n - \eta)} \; := \; \Omega^\infty 
\C P^\infty_{-1}\;.
\]

\noindent
{\bf 2.2} Madsen and Tillmann show their construction factors through 
an infinite loopspace map
\[
\Z \times B\Gamma^+_\infty \to \Omega^\infty \C P^\infty_{-1} \to 
Q(\C P^\infty_+) \;,
\]
in which the last arrow is defined by collapsing the bottom two-cell in a 
cofibration 
\[
S^{-2} \to \C P^\infty_{-1} \to \C P^\infty_+ \;.
\]
The fiber $\Omega^2 QS^0$ of the induced map of loop spaces is torsion, 
so the rational cohomology of $\Omega^\infty \C P^\infty_{-1}$ is isomorphic 
to the algebra of symmetric functions on the subspace of non-negative powers
in $t^*_\T H\Q$. This algebra is thus canonically isomorphic to the Fock 
representation [10 \S 2.2] of the Heisenberg algebra of that symplectic 
module; but this representation possesses a canonical Virasoro action, 
defining a homomorphism 
\[
H_*(\Z \times B\Gamma^+_\infty,\Q) \to {\rm Symm}(H^*(\C P^\infty_+)) \in 
(\GG^-_{1/2}-{\rm representations}) \;.
\]
In Kontsevich-Witten theory the usual generators $b_k \in H_*(\C P^\infty_+), 
k \geq 0$, map to symmetric functions
\[
\tr \; \gamma_{-k-\half}(\Lambda^2) \sim -(2k-1)!! \; \tr \; \Lambda^{-2k-1} 
= t_k(\Lambda)
\]
of a positive-definite Hermitian matrix $\Lambda$; this leads to a 
construction of the appropriate twisted Virasoro representation in terms
of Schur $Q$-functions [2,6]. \bigskip

\noindent
{\bf 2.3} The homomorphism
\[
\lim MU^{*+n-2}(Th(n - \eta)) \to MU^{*-2}(B\Diff(F)) 
\]
defined on cobordism by the Madsen-Tillmann construction sends the Thom 
class to a kind of Euler class: according to Quillen's conventions, the 
Thom class is the zero-section of $(n - \eta)$, regarded as a cobordism 
class of maps between manifolds. Its image is the class defined by
the fiber product
$$\xymatrix{
{Z_n} \ar[r] \ar[d]& {\gr_{2,n}} \ar[d]  \\
{\R^n_+ \wedge_{\Diff} \Emb(F)} \ar[r]& {\gr_{2,n}^{(n - \eta)}} \;;}
$$ 
it is the space of equivalence classes, under the action of
$\Diff(F)$, of pairs $(x,\phi)$, with $x \in \phi(F) \subset \R^n$
a point of the surface (ie, in the zero-section of $\nu$), and $\phi$ 
an embedding. The image is thus the element
\[
[Z_n \to \R^n_+ \wedge_{\Diff} \Emb] \mapsto MU^{n-2}(S^n B\Diff(F)) 
\]
defined by the tautological family $F \times_{\Diff} E\Diff(F)$ of surfaces 
over the classifying space of the diffeomorphism group. This class is 
primitive in the Hopf-like structure defined by gluing [9 \S 2.2], so 
the class
\[
\Phi = \exp (th(-\eta)v) \in MU^0_\Q(\Omega^\infty \C P^\infty_{-1})[[v]]
\]
of finite unordered configurations of points on the universal surface 
(with $v$ a book-keeping indeterminate) defines a multiplicative 
transformation
\[
\tilde \Phi_* : H_*(Q(\C P^\infty_+),\Q) \to H_*(MU,\Q[[v]])
\]
with properties analogous to the Chern character of a vector bundle. It
sends the Fock representation described above to an algebra of cohomological 
characteristic numbers. \bigskip

\noindent
{\bf 2.4} Kontsevich-Witten theory uses a more sophisticated 
configuration space, which maps the rational homology of $Q(\coprod_{g 
\geq 0} \M_g)$ (suitably interpreted, for small $g$) to a similar ring of 
characteristic numbers; this homology contains a fundamental class
\[
[Q(\coprod \M_g)] = \exp(\sum_{g \geq 0} [\M_g]v^{3(g-1)})
\]
for the moduli space of not necessarily connected curves. Witten's 
tau-function is the image of this `highest-weight' vector under the analog 
of $\tilde \Phi_*$; it is invariant under the subalgebra of Virasoro
generated by operators $L_k$ (of cohomological degree $k$) with $k \geq -1$.

\section{Afterthoughts, and possible generalizations}
\bigskip

\noindent
{\bf 3.1} Witten has proposed a generalization of 2D topological
gravity which encompasses surfaces with higher spin structures:
for a closed smooth surface $F$ an $r$-spin structure is roughly
a complex line bundle $L$ together with a fixed isomorphism
$L^{\otimes r} \cong T_F$ of two-plane bundles, but for surfaces
with nodes or marked points the necessary technicalities are
formidable [5]. The group of automorphisms of such a structure is
an extension of its group of diffeomorphisms by the group of
$r$th roots of unity, and there is a natural analog of the group
completion of the category defined by such surfaces. The
generalized Madsen-Tillmann construction maps this loopspace to
the Thom spectrum $Th(-\eta^r)$, and it is reasonable to expect
that this map is equivariant with respect to automorphisms of the
group of roots of unity. This fits with some classical homotopy 
theory: if (for simplicity) $r=p$ is prime, multiplication by an 
integer $u$ relatively prime to $p$ in the $H$-space structure of 
$\C P^\infty$ defines a morphism
\[
Th(-\eta^p) \to Th(-\eta^{up})
\]
of spectra, and the classification of fiber-homotopy equivalences
of vector bundles yields an equivalence of $Th(-\eta^{up})$ with
$Th(-\eta^p)$ after $p$-completion. There is an analogous
decomposition of $t_TH\Z_p$ and a corresponding decomposition of
the associated Fock representations [10 \S 2.4]. \bigskip

\noindent
{\bf 3.2} In an extension of her work Tillmann also considers categories 
of surfaces mapped to some parameter space $X$, which has interesting 
connections with both Tate and quantum cohomology. When $X$ is a smooth 
compact almost-complex manifold, its Hodge-deRham cohomology admits a 
natural action of the Lie algebra ${\bf sl}_2(\R)$, generated [12 IV \S 4] 
by the Hodge dimension operator $H$, multiplication by the first Chern 
class $c_1(X) = E$, and its adjoint $F = *E*$. Recently Givental [3 \S 8.1] 
has shown that earlier work of (the schools of) Eguchi, Dubrovin, and others 
on the Virasoro structure of quantum cohomology can be formulated in terms 
of (what I like to think of as) $t^*_\T H(X,\R)$, polarized by the twisted 
involution
\[
I_{Giv} = \e^{\half H} \e^{-E} I \e^E \e^{-\half H}  \;.
\]
It would be very interesting if this polarization could be understood
in terms of the equivariant geometry of the free loopspace of $X$. 
\bigskip

\noindent
{\bf 3.3} I owe thanks to R. Cohen, E. Getzler, A. Givental, J. Greenlees,
I. Madsen, N. Strickland, and U. Tillmann (at least), for their forebearance 
in the face of my continued misunderstanding of things they have tried 
patiently to tell me. I hope I'm finally starting to get it right. 
\bigskip

\bibliographystyle{amsplain}

\end{document}